\input amssym.def
\input amssym
\magnification=1200
\parindent0pt
\hsize=16 true cm
\baselineskip=13  pt plus .2pt
$ $

\def\R{\Bbb R}
\def\C{\Bbb C}

\centerline {\bf On Jordan type bounds for finite groups of diffeomorphisms}
\centerline {\bf of 3-manifolds and Euclidean spaces}

\bigskip  \bigskip

\centerline {Bruno P. Zimmermann}

\medskip

\centerline {Universit\`a degli Studi di Trieste}

\centerline {Dipartimento di Matematica e Geoscienze}

\centerline {34127 Trieste, Italy}

\bigskip  \bigskip

Abstract.  {\sl  By a classical result of Jordan,  each finite subgroup $G$ of
${\rm GL}_n(\C)$ has an abelian subgroup whose index in $G$ is bounded by a
constant depending only on $n$. We consider the problem if this remains true
for finite subgroups $G$ of the diffeomorphism group of a smooth manifold, and
show that it is true for all compact 3-manifolds as well as for Euclidean
spaces $\R^n$, $n \le 6$. The question remains open at present e.g. for
odd-dimensional spheres $S^n$, $n \ge 5$ and for Euclidean spaces $\R^n$, $n
\ge 7$.}

\bigskip \bigskip

{\bf 1. Introduction}

\medskip

By a classical result of Jordan,  each finite subgroup $G$ of ${\rm GL}_n(\C)$
has an abelian subgroup $A$ whose index in $G$ is bounded by a constant
depending only on $n$ (see [C] for the optimal bound for each $n$).  Recently
there has been much interest in generalizations, replacing ${\rm GL}_n(\C)$ by
more general geometrically interesting groups such as diffeomorphism groups of
smooth manifolds ([MR1,2],[P1]), automorphism groups of algebraic varieties and
the Cremona groups of birational self-maps of the affine $n$-dimensional space
(cf. [P2], [Se, Theoreme 3.1]).

\medskip

Following [P1,2], we say that a group $E$ is a {\it Jordan group} or has the
{\it Jordan property} if there exits a constant such that every finite subgroup
$G$ of $E$ has an abelian  subgroup of index bounded by this constant. For a
smooth manifold $M$, let ${\rm Diff}(M)$ denote its diffeomorphism group. The
present paper is motivated by the following general:

\medskip

Question: For which (classes of) smooth manifolds $M$ is ${\rm Diff}(M)$ a
Jordan group?

\medskip

Whereas this is in general not true for non-compact manifolds ([P1]), it has
been conjectured that it is true for compact manifolds (see [MR1,2]); however,
it should be true e.g. also for ${\rm Diff}(\R^n)$.

\medskip

Note that a Jordan group contains only finitely many finite non-abelian simple
subgroups, up to isomorphism; in this regard, it has been shown in [GZ] that
${\rm Diff}(S^n)$ contains only finitely many finite non-abelian simple
subgroups, up to isomorphism  (and, more generally, for any closed homology
$n$-sphere, see also [Z1]). It has been shown in [MR1] that ${\rm Diff}(M)$ is
a Jordan group if $M$ is a compact manifold without odd cohomology; in
particular, ${\rm Diff}(S^n)$ is a Jordan group for even dimensions $n$, but
this remains open for odd dimensions $n \ge 5$.

\medskip

On the basis of the geometrization of 3-manifolds and results of Kojima [K] and
the author [Z2], in our first main result we consider the case of compact
3-manifolds:

\bigskip

{\bf Theorem 1.}  {\sl  ${\rm Diff}(M)$ is a Jordan group for compact 3-manifolds $M$.}

\bigskip

In dimension three, this leaves open the question for non-compact 3-manifolds.
Concerning dimension four, it is shown in [P2] that there are noncompact,
simply connected,  smooth 4-manifolds $M$ such that ${\rm Diff}(M)$ is not a
Jordan group. On the other hand, it is shown in [MR2] that, for compact smooth
4-manifolds $M$ with non-zero Euler characteristic, ${\rm Diff}(M)$ is a Jordan
group. The case of the 4-sphere $S^4$ is considered in [MeZ1,2] where it is
shown that, up to 2-fold extension in the case of solvable groups, any finite
group with an orientation-preserving smooth action on $S^4$ (or on any homology
4-sphere) is isomorphic to a subgroup of the orthogonal group SO(5), presenting
also a short list of such groups.

\medskip

Next we consider Euclidean spaces $\R^n$. The following is proved in [GMZ]:

\bigskip

{\bf Theorem 2.}([GMZ])  {\sl Let $G$ be a finite subgroup of ${\rm Diff}(\R^n)$ (or of
${\rm Diff}(M)$, for any acyclic $n$-manifold $M$).  Suppose that $n \le 4$; then $G$ is
isomorphic to a subgroup of the orthogonal group
${\rm O}(n)$. In particular, the classical Jordan bound applies to $G$, so
${\rm Diff}(\R^n)$ is a Jordan group for $n \le 4$.}

\bigskip

In [GMZ] the case of finite groups of diffeomorphisms of $\R^5$ is also considered; the
classification in this case is not complete but the results imply easily that also ${\rm
Diff}(\R^5)$ is a Jordan group (more generally, the results in [GMZ] apply to arbitrary
acyclic manifolds).

\medskip

A main tool for the proof of our second main result is a recent
group-theoretical result of Mundet i Riera and Turull [MT] (on the basis of the
classification of the finite simple groups).

\bigskip

{\bf Theorem 3.}  {\sl ${\rm Diff}(\R^5)$ and  ${\rm Diff}(\R^6)$ are Jordan
groups (and, more generally, ${\rm Diff}(M)$ for any acyclic 5- or 6-manifold
$M$).}

\bigskip

We will present the proof of Theorem 3 for the new case $n = 6$; the same proof
works also for $n = 5$ where it is, in fact, easier. As noted above, the proof
for $n=6$ uses the full classification of the finite simple groups; the proof
for $n=5$ instead requires "only" a smaller part of the classification of the
finite simple groups (the Gorenstein-Harada classification of the finite simple
groups of sectional 2-rank at most four), see [GMZ], [Z1].

\medskip

Two interesting cases where the Jordan property is not known at present are
those of  ${\rm Diff}(S^5)$ and ${\rm Diff}(\R^7)$. However, it seems likely
that ${\rm Diff}(S^n)$ and ${\rm Diff}(\R^n)$ are Jordan groups for all values
of $n$.

\bigskip  \bigskip

{\bf 2. Proof of Theorem 1}

\medskip

It is easy to see that, if  $\tilde M$ is a finite covering of $M$ such that
${\rm Diff}(\tilde M)$ is a Jordan group then also ${\rm Diff}(M)$ is a Jordan
group. So it is sufficient to consider the case of orientable manifolds, and
also of orientation-preserving finite group actions $G$ (passing eventually to
a subgroup of index two of $G$). Also, it is sufficient to consider the case of
closed manifolds since, for a compact manifold $M$ with non-empty boundary, one
can reduce to the closed case by taking the double of $M$ along the boundary,
 doubling also a given finite group action on $M$.

\medskip

So let $M$ be a closed orientable 3-manifold and $G$ a finite group of
orientation-preserving diffeomorphisms of $M$. If $\pi_1(M)$ is finite then, by
the geometrization of 3-manifolds after Perelman, $M$ is a spherical 3-manifold
and finitely covered by $S^3$; also, any finite group of diffeomorphisms of $M$
is conjugate to a linear (orthogonal) action. By the classical Jordan bound for
linear groups, ${\rm Diff}(S^3)$ is a Jordan group, and hence also ${\rm
Diff}(M)$ is a Jordan group.

\medskip

Assume next that $M$ is irreducible and has infinite fundamental group; again by the
geometrization of 3-manifolds, we can assume that $M$ is a geometric. Then, if $M$
does not admit a circle action, by [K, Theorem 4.1] there is a bound on the order of
finite subgroups of ${\rm Diff}(M)$ and we are done.

\medskip

Suppose that $M$ has a circle action and infinite fundamental group. Then $M$
is a Seifert fiber space, and by the geometrization of finite group actions on
Seifert fiber spaces ([MS]), we can assume that the action of the finite group
$G$ of diffeomorphisms of $M$ is geometric, and in particular fiber-preserving
and normalizing the $S^1$-action of $M$. Considering a suitable finite covering
of $M$, we can moreover assume that $M$ has no exceptional fibers, and hence
that the base space of the Seifert fibration (the quotient of the $S^1$-action)
is a closed orientable surface $B$ without cone points.  The finite group $G$
projects to a finite group $\bar G$ of diffeomorphisms of the base-surface $B$,
and we can again assume that $\bar G$ is orientation-preserving.

\medskip

If $B$ is a hyperbolic surface (of genus $g \ge 2$) then, by the formula of
Riemann-Hurwitz, the order of the finite group $\bar G$ of diffeomorphisms of
$B$ is bounded, and hence $G$ has a finite cyclic subgroup of bounded index
(the intersection of $G$ with the $S^1$-action).

\medskip

If $B$ is a torus $T^2$ then there are two cases. First, $M$ may be a
3-dimensional torus $T^3$; this acts by rotations on itself. Since the action
of $G$ is geometric, the subgroup $G_0$ of $G$ acting trivially on the
fundamental group is a subgroup of the $T^3$-action and hence abelian of rank
at most three (see [Sc] for the geometries of 3-manifolds and their isometry
groups). The factor group $G/G_0$ acts faithfully on the fundamental group
$\Bbb Z^3$ of the 3-torus and is  isomorphic to a subgroup of ${\rm GL}_3(\Bbb
Z)$. Since, by a well-known result of Minkowski, there is a bound on the finite
subgroups of ${\rm GL}_n(\Bbb Z)$ for each $n$, $G$ has an abelian subgroup
$G_0$ of bounded index.

\medskip

If $M$ fibers over $T^2$ but is not a 3-torus then it belongs to the nilpotent
geometry Nil given by the Heisenberg group (see again [Sc]). Now the subgroup
$G_0$ of $G$ acting trivially on the fundamental group, up to inner
automorphisms, is a cyclic subgroup of the $S^1$ action on $M$, and $G/G_0$
injects into the outer automorphism group ${\rm Out}(\pi_1M)$ of the
fundamental group. The fundamental group of $M$ has a presentation
$$\pi_1M = \; <a,b,t \mid  [a,b]=t^k, \;\; [a,t]=[b,t]=1>,$$
with $k \ne 0$. Now an easy calculation shows that the subgroup of the outer
automorphism group of $\pi_1M$ inducing the identity of the factor group
$\pi_1M)/<t> \; \cong \; \Bbb Z^2$ is finite. Since the orders of finite
subgroups of ${\rm GL}_2(\Bbb Z)$ are also bounded, $G$ has a finite cyclic
subgroup $G_0$ of bounded index.

\medskip

Finally, if the base-surface is the 2-sphere then either $M$ has finite
fundamental group and is a spherical manifold, or homeomorphic to $S^2 \times
S^1$ (and hence non-irreducible). We note that $S^2 \times S^1$ belongs to the
$(S^2 \times \Bbb R$)-geometry, one of Thurston's eight 3-dimensional
geometries; this is the easiest of the eight geometries and can be easily
handled, see [Sc] for the isometry group of this geometry.

\medskip

Summarizing, we have shown that for any closed irreducible 3-manifold $M$ (and
also for $S^2 \times S^1$), ${\rm Diff}(M)$ is a Jordan group.

\medskip

Suppose that $M$ is non-irreducible but not $S^2 \times S^1$. If $M$ has a
summand other than lens spaces and $S^2 \times S^1$ then, by [K, Theorem 4.2],
the orders  of finite diffeomorphism groups of $M$ are again bounded and we are
done.

\medskip

Suppose next that $M$ is a connected sum $\sharp_g (S^2 \times S^1)$ of $g$
copies of $S^2 \times S^1$, with $g > 1$. By [Z2], $G$ has a finite cyclic
normal subgroup (the subgroup acting trivially on the fundamental group, up to
inner automorphisms) such that the order of the factor group is bounded by a
polynomial which is quadratic in $g$, so we are done also in this case.
Finally, if $M$ is a connected sum of lens spaces, including $S^2 \times S^1$,
then $M$ has a finite covering by a 3-manifold of type $\tilde M = \sharp_g
(S^2 \times S^1)$ as before. Now ${\rm Diff}(\tilde M)$ is a Jordan group and
hence also ${\rm Diff}(M)$.

\medskip

We have considered all possibilities for $M$ and completed the proof of Theorem 1.

\bigskip  \bigskip
\vfill \eject

{\bf 3. Proof of Theorem 3}

\medskip

We prove the theorem for $n=6$; for  $n = 5$ the theorem follows from [GMZ,
Theorem 3], and also a shorter version of the following proof for $n = 6$
applies.

\medskip

We want to show  that ${\rm Diff}(\R^6)$ is a Jordan group, i.e. that there is
a constant such that every finite subgroup $G$ of ${\rm Diff}(\R^6)$ has an
abelian subgroup of index bounded by this constant. By the main result of [MT],
if this is true for all finite subgroups $G$ of ${\rm Diff}(\Bbb R^6)$  which
are a semidirect product $G = P \rtimes Q$, for a finite normal $p$-group $P$
and a finite $q$-group $Q$, with distinct primes $p$ and $q$,  then it is true
for all finite subgroups $G$ of ${\rm Diff}(\Bbb R^6)$ (this uses the
classification of the finite simple groups). So we have to consider only groups
of type  $G = P \rtimes Q$: given such a group, we have to find an abelian
subgroup $A$ of $G$ whose index is bounded by a constant not depending on the
specific group.

\medskip

Let $G = P \rtimes Q$ be as before; we can assume that the action of $G$ is
orientation-preserving.  By general Smith fixed point theory, a finite
$q$-group acting on $\Bbb R^n$ (or on any acyclic $n$-manifold) has non-empty
fixed point set (see [B], [GMZ, section 2]). So $Q$ has a global fixed point
and is isomorphic to a subgroup of the orthogonal group SO(6) (considering the
induced linear action on the tangent space of a global fixed point). Hence, by
the classical Jordan bound for linear groups, we may assume that $Q$ is an
abelian $q$-group.

\medskip

Let $F$ denote the fixed point set of $P$; since $P$ is normal, $F$ is invariant under
the action of $Q$ and, since the action is orientation-preserving, $F$ is a submanifold
of dimension at most four (i.e., of codimension at least two).

\medskip

Suppose first that $F$ has dimension four.  Then $P$ acts as a group of
rotations around its fixed point set $F$ and hence is a cyclic group
(isomorphic to a subgroup of SO(2)). By conjugation, every element of $G$ acts
as $\pm$-identity on $P$ (conjugates a minimal rotation to a minimal rotation).
Let $G_0$ be the subgroup of index one or two of $G$ acting trivially on $P$,
and let $Q_0$ be its image in $Q$. Then $G_0 \cong P \times Q_0$ is an abelian
subgroup of index at most two in $G$, so we are done in this case.

\medskip

Now suppose that the fixed point set $F$ of $P$ has dimension three (and also
codimension three). This implies that $p=2$ since, if $p$ is odd, by an
inductive argument on the $p$-group $P$, its fixed point set $F$ has even
codimension. Considering the action of $P$ on a 3-ball transverse to $F$ in
some point, $P$ is a subgroup of the orthogonal group SO(3) and hence
isomorphic to a cyclic or dihedral 2-group.

\medskip

If $P \cong \Bbb Z_2 \times \Bbb Z_2$ is isomorphic to the Klein 4-group then
the subgroup $G_0$ of $G$ acting by conjugation trivially on $P$ has index at
most three in $G$ (since $Q$ is a $q$-group of odd order) and is abelian, so we
are done. If $P$ is a cyclic 2-group then its automorphism group is also a
2-group; since $Q$ has odd order, $G$ acts by conjugation trivally on $P$, so
$G$ is abelian and we are done. If $P$ is a dihedral $2$-group of order at
least eight then it has a cyclic characteristic subgroup $P_0$ of index two, so
$G$ has a subgroup $G_0 = P_0 \rtimes Q$ of index two; by the previous case,
$G_0$ is abelian and we are done.

\medskip

Suppose next that $F$ has dimension two. By Smith fixed point theory, $F$ is an
acyclic manifold mod $p$ (i.e., for homology with coefficients in $\Bbb Z_p$).
Since $F$ has dimension two, it is in fact acyclic also for integer
coefficients (see [GMZ, proof of Lemma 3]). Then the finite $q$-group $Q$ has a
fixed point in $F$, and hence $G$ has a global fixed point. Now $G$ is
isomorphic to a subgroup of SO(6), so we are done by the classical Jordan
bound.

\medskip

The cases that $F$ has dimension one or zero are similar.

\medskip

This completes the proof of Theorem 3.

\bigskip

Remark.  Considering the next case of $\Bbb R^7$, if the fixed point set $F$ of
$P$ has codimension two or three, or if it has dimension at most two, the proof
works exactly as before. The case we cannot handle at present is when $F$ has
dimension three (and codimension four). In this case $P$ is isomorphic to a
subgroup of SO(4), e.g. isomorphic to $\Bbb Z_p$ or $\Bbb Z_p \times \Bbb Z_p$,
and we don't know how to bound the index of the subgroup of $G$ (or $Q$) acting
trivially on $P$  (independent of the prime $p$).

\bigskip \bigskip

\centerline {\bf References}

\bigskip

\item {[B]} G. Bredon, {\it Introduction to Compact Transformation Groups,}
Academic Press, New York 1972

\item {[C]}  M.J. Collins, {\it  On Jordan's theorem for complex linear groups,}
J. Group Theory 10,  411-423 (2007)

\item {[GMZ]} A. Guazzi, M. Mecchia, B. Zimmermann, {\it On finite groups acting on
acyclic low-dimensional manifolds,}  Fund. Math. 215,  203-217  (2011)

\item {[GZ]} A. Guazzi, B. Zimmermann, {\it On finite simple groups acting on
homology spheres,}  Mo-natsh. Math. 169,  371-381 (2013)

\item {[K]} S. Kojima, {\it Bounding finite groups acting on 3-manifolds,}
Math. Proc. Camb. Phil. Soc. 96,  269-281  (1984)

\item {[MeZ1]} M. Mecchia, B. Zimmermann, {\it On finite simple and nonsolvable groups
acting on homology 4-spheres,} Top. Appl. 153,  2933-2942  (2006)

\item {[MeZ2]} M. Mecchia, B. Zimmermann, {\it On finite  groups acting on homology
4-spheres and finite subgroups of ${\rm SO}(5)$,}  Top. Appl. 158,  741-747  (2011)

\item {[MS]}  W.H. Meeks, P. Scott, {\it  Finite group actions on 3-manifolds,}
Invent. math. 86, 287-346  (1986)

\item {[MR1]} I. Mundet i Riera, {\it Finite groups acting on manifolds without odd
cohomology,}  arXiv: 1310.6565

\item {[MR2]} I. Mundet i Riera, {\it Finite group actions on 4-manifolds with
nonzero Euler characteristic,}
arXiv:1312.3149

\item {[MT]} I. Mundet i Riera, A. Turull, {\it Boosting an analogue of Jordan's
theorem for finite groups,}
arXiv:1310.6518

\item {[P1]} V.L. Popov, {\it Finite subgroups of diffeomorphism groups,}
arXiv:1310.6548

\item {[P2]} V.L. Popov, {\it Jordan groups and automorphism groups of algebraic
varieties,}  arXiv: 13007.5522

\item {[Sc]}  P. Scott, {\it The geometries of 3-manifolds,}
Bull. London Math. Soc. 15, 401-487  (1983)

\item {[Se]} J.-P. Serre, {\it Le groupe de Cremona et ses sous-groupes finis,}
Sem. Bourbaki 1000, 75-100 (2008)

\item {[Z1]} B. Zimmermann, {\it On finite groups acting on spheres and finite
subgroups of orthogonal groups,} Sib. Electron. Math. Rep. 9, 1 - 12  (2012)
(http://semr.math.nsc.ru)

\item {[Z2]}   B. Zimmermann,  {\it  On finite groups acting on a connected sum of
3-manifolds $S^2 \times S^1$,}  arXiv:1202.5427  (to appear in Fund. Math.
2014)

\bye